\documentclass[12pt,english]{article}
\usepackage{amsmath,amssymb,euscript,amsthm,amsfonts}
\usepackage[cp1251]{inputenc}
\usepackage[english,ukrainian]{babel}
\begin{document}

\begin{center}
\textbf{\Large Order estimates of the uniform approximations by
Zygmund sums on the classes of convolutions of periodic functions}
\end{center}
\vskip0.5cm

\begin{center}

Serdyuk A.S., Hrabova U.Z. \\ \emph{\small Institute of Mathematics,
Kiev \\ Volyn national university of Lesya Ukrainka, Lutsk}
\end{center}
\vskip0.5cm

\begin{center}

 {\bf Abstract}
\end{center}
\vskip0.3cm

We establish the exact-order estimates of uniform approximations by
the Zygmund sums $Z^{s}_{n-1}$ (that is trigonometric polynomials of
the form
\linebreak$Z^{s}_{n-1}(f;t):=\frac{a_{0}}{2}+\sum_{k=1}^{n-1}\left(1-\left(\frac{k}{n}\right)^{s}\right)\times
(a_{k}(f)\cos kt+b_{k}(f)\sin kt), s>0,$ where $a_{k}(f)$ and
$b_{k}(f)$ are the Fourier coefficients of $f\in L_{1}$) of
$2\pi$-periodic continuous functions $f$ from the classes
$C^{\psi}_{\beta,p}$. These classes are defined by the convolutions
of functions from the unit ball in the space
  $L_{p}$, $1\leq p<\infty$, with generating fixed kernels
 $\Psi_{\beta}(t)=\sum_{k=1}^{\infty}\psi(k)\cos\left(kt+\frac{\beta\pi}{2}\right)$,
  $\Psi_{\beta}\in L_{p'}$, $\beta\in \mathbb{R}$, $1/p+1/p'=1$. We additionally assume that the product
   $\psi(k)k^{s+1/p}$ is generally monotonically increasing with the rate of some power  function, and, besides, for $1< p<\infty$  it holds that
     $\sum_{k=n}^{\infty}\psi^{p'}(k)k^{p'-2}<\infty$,
 and for $p=1$ the following condition is true
 $\sum_{k=n}^{\infty}\psi(k)<\infty$.

  It is
shown that under these conditions Zygmund sums $Z^{s}_{n-1}$ and
Fejer sums \linebreak$\sigma_{n-1}=Z^{1}_{n-1}$ realize the order of
the best uniform approximations by trigonometric polynomials of
these classes, namely for $1<p<\infty$
$${E}_{n}(C^{\psi}_{\beta,p})_{C}\asymp{\cal
E}\left(C^{\psi}_{\beta,p};
Z_{n-1}^{s}\right)_{C}\asymp\Big(\sum_{k=n}^{\infty}\psi^{p'}(k)k^{p'-2}\Big)^{1/p'},
\ \frac{1}{p}+\frac{1}{p'}=1,$$ and for $p=1$
$$
{E}_{n}(C^{\psi}_{\beta,1})_{C}\asymp{\cal
E}\left(C^{\psi}_{\beta,1}; Z_{n-1}^{s}\right)_{C}\asymp
{\left\{{\begin{array}{l l}
\sum\limits_{k=n}^{\infty}\psi(k), & \cos \frac{\beta\pi}{2}\neq 0; \\
\psi(n)n, &\cos \frac{\beta\pi}{2}= 0,
\end{array}} \right.}
$$
where
 $${E}_{n}(C^{\psi}_{\beta,p})_{C}:=\sup_{f\in
C^{\psi}_{\beta,p}}\inf\limits_{t_{n-1}\in\mathcal{T}_{2n-1}}\|f-t_{n-1}\|_{C},
$$
and $\mathcal{T}_{2n-1}$ is the subspace of trigonometric
polynomials $t_{n-1}$ of order  $n-1$ with real coefficients,
$${\cal E}\left(C^{\psi}_{\beta,p};
Z_{n-1}^{s}\right)_{C}:=\mathop{\sup}\limits_{f\in
C^{\psi}_{\beta,p}}\|f-Z^{s}_{n-1}(t)\|_{C}.$$


\vskip0.6cm

\vskip1cm

 \normalsize \vskip 3mm

 \section{Notations, definitions and   auxiliary  statements }

Denote by $L_{p}$, $1\leq p\leq\infty$,
 the space of  $2\pi$--periodic summable  on $[0,2\pi]$ functions
$f$ with the norm
$$\|f\|_{p}={\left\{
{\begin{array}{l
l}\Big(\int\limits_{0}^{2\pi}|f(t)|^{p}dt\Big)^{1/p}, &
1\leq p<\infty; \\
\mathop{\rm{ess}\sup}\limits_{t}|f(t)|,  & p=\infty,
\end{array}} \right.}
$$
and by  $C$ the space of $2\pi$--continuous periodic functions in
which the norm is defined by equality
$\|f\|_{C}=\max\limits_{t}|f(t)|$.

Let  $f\in L_{1}$ and
$$
S[f](x)=\frac{a_0}{2}+\sum_{k=1}^{\infty}(a_k(f) \cos kx + b_k
(f)\sin kx),
$$
be  the Fourier series of function $f$.

If  for the  sequence  $\psi(k)\in \mathbb{R}$ and fixed number
$\beta\in \mathbb{R}$
 the series
$$
\sum_{k=1}^{\infty}\frac{1}{\psi\left(k\right)}\bigg(a_{k}(f)\cos\Big(kx+\frac{\beta\pi}{2}\Big)
+b_{k}(f)\sin\Big(kx+\frac{\beta\pi}{2}\Big)\bigg)
$$
 is the Fourier series of a
summable function $\varphi$, then this function is called as
$(\psi,\beta)$-derivative of the  function $f(x)$  and is denoted by
$f_{\beta}^{\psi}(x)$. A set of functions $f(x)$, for which this
condition is satisfied is denoted by $L_{\beta}^{\psi}$, and subset
all continuous functions from   $L^{\psi}_{\beta}$ is denoted by
$C^{\psi}_{\beta}$.

If $f\in L^{\psi}_{\beta}$ and furthermore $f^{\psi}_{\beta}\in
\mathfrak{N}$, where $\mathfrak{N}\subset L_{1}$, then we write that
$f\in L^{\psi}_{\beta}\mathfrak{N}$. Let us put
$L^{\psi}_{\beta}\mathfrak{N}\cap C= C^{\psi}_{\beta}\mathfrak{N}$.
The concept of  $(\psi,\beta)$-derivative is a natural
generalization of the concept of   $(r,\beta)$-derivative in the
Weyl–Nagy sense  and coincides  almost everywhere with the last
one, when
 $\psi(k)=k^{-r}$, $ r>0$, namely, if  $\psi(k)=k^{-r}$, $ r>0$, then  ${L}_{\beta}^{\psi}\mathfrak N=W^{r}_{\beta}{\mathfrak N}$, and,
 $f^{\psi}_{\beta}=f^{r}_{\beta}$,
   where $f^{r}_{\beta}$ is  the derivative in the Weyl–Nagy sense,
   and
  $W^{r}_{\beta}{\mathfrak N}$ are  the Weyl-Nagy classes  \cite{Nady}, \cite{S1}. In the case, when $\beta=r$, the classes $W^{r}_{\beta}{\mathfrak N}$
  are the well known   Weyl classes
    $W^{r}_{r}{\mathfrak N}$, while  the derivatives $f^{r}_{\beta}$ coincide almost everywhere  with the derivatives in the sense of Weyl
   $f^{r}_{r}$.
If, in addition,  $\beta=r$, $r\in\mathbb{N}$, then $f_{\beta}^{r}$
 coincide almost everywhere with the usual derivatives  $f^{(r)}$ of
the order
   $r$ of the function $f$ ($f^{r}_{\beta}=f^{r}_{r}=f^{(r)}$) and at the same time  $W^{r}_{\beta}{\mathfrak N}=W^{r}_{r}{\mathfrak N}=W^{r}{\mathfrak N}$.

According to the Statement 3.8.3  from    \cite{S1}, if the series
\begin{equation}\label{rff}
\sum_{k=1}^{\infty}\psi(k)\cos \big(kt-\frac{\beta\pi}{2}\big), \
\beta\in
    \mathbb{R}
\end{equation}
is the Fourier series of the function  $\Psi_{\beta}\in L_{1}$,
 then the elements  $f$ of the classes  $L_{\beta}^{\psi}\mathfrak N$ for almost
 every
$x\in\mathbb{R}$   are represented as a convolution
\begin{equation}\label{zgo}
f(x)=\frac{a_{0}}{2}+(\Psi_{\beta}\ast\varphi)(x)=\frac{a_{0}}{2}+
\frac{1}{\pi}\int\limits_{-\pi}^{\pi}\Psi_{\beta}(x-t)\varphi(t)dt,
\ a_{0}\in\mathbb{R}, \varphi\perp1, \  \varphi\in\mathfrak N,
\end{equation}
where   $\varphi$  almost everywhere coincides with
$f^{\psi}_{\beta}$.

As sets  $\mathfrak N$ we will consider the unit balls of the spaces
$L_{p}$:
 $$
U_{p}=\{\varphi\in L_{p}:\|\varphi\|_{p}\leq~1 \}, \ 1\leq
p\leq\infty.
$$
Then put:
 $L^{\psi}_{\beta,p}:=L^{\psi}_{\beta}U_{p}$, \
 $C^{\psi}_{\beta,p}:=C^{\psi}_{\beta}U_{p}$, $W^{r}_{\beta,p}:=W^{r}_{\beta}U_{p}$.

According to the Statement  1.2, from  \cite{S1}, if the fixed
kernel $\Psi_{\beta}$ of the classes  $L^{\psi}_{\beta,p}$ and
$C^{\psi}_{\beta,p}$ satisfies the inclusion $\Psi_{\beta}\in
L_{p'}$, $\frac{1}{p}+\frac{1}{p'}=1$, $1\leq p\leq\infty$, then the
convolutions of the form (\ref{zgo}) are continuous functions, where
 $\mathfrak N=U_{p}$. It is clear that in this case for   $f\in C^{\psi}_{\beta,p}$ the equality
(\ref{zgo}) is fulfilled for all $x\in\mathbb{R}$.

We assume that the sequences  $\psi(k)$ are traces on the set of
natural numbers  $\mathbb{N}$ of some positive continuous convex
downwards functions $\psi(t)$ of the continuous argument $t\geq 1$,
 that tends  to zero for
$t\rightarrow\infty$. The set of all such functions  $\psi(t)$ is
denoted by $\mathfrak{M}$.

To classify functions $\psi$ from $\mathfrak{M}$ on their speed of
decreasing to zero it is convenient to use the  following
characteristic :
\begin{equation}\label{alfa}
\alpha(t)=\alpha(\psi;t)=\frac{\psi(t)}{t|\psi'(t)|}, \
\psi'(t):=\psi'(t+0).
\end{equation}
With its help we consider  the following subsets  of the set
$\mathfrak{M}$  (see, e.g., \cite{S1})
$$
\mathfrak{M}_0:=\{\psi\in\mathfrak{M}: \ \ \exists K>0 \ \ \forall
t\geq1 \ \ 0<K\leq\alpha(\psi;t)\},
$$
$$
\mathfrak{M}_C:=\{\psi\in\mathfrak{M}: \ \ \exists K_{1}, K_{2}>0 \
\ \forall t\geq1 \ \ 0<K_{1}\leq\alpha(\psi;t)\leq K_{2}\}.
$$
It is clear that $\mathfrak{M}_{C}\subset\mathfrak{M}_{0}$.

 Zygmund sums of the order $n-1$
 of the function $f\in L_{1}$ are the trigonometric polynomials of the form
\begin{equation}\label{sz}
Z_{n-1}^{s}(f;t)=\frac{a_{0}}{2}+\sum_{k=1}^{n-1}\left(1-\left(\frac{k}{n}\right)^{s}\right)(a_{k}(f)\cos
kt+b_{k}(f)\sin kt), s>0,
\end{equation}
where $a_{k}(f)$ and $b_{k}(f)$ are Fourier coefficients of the
function $f$.

In the case $s=1$ polynomials  $Z^{s}_{n-1}$
 are  Fejer sums:
$Z^{1}_{n-1}=\sigma_{n-1}$
\begin{equation}\label{sf}
\sigma_{n-1}(f;t)=\frac{a_{0}}{2}+\sum_{k=1}^{n-1}\left(1-\frac{k}{n}\right)(a_{k}(f)\cos
kt+b_{k}(f)\sin kt).
\end{equation}

In this paper we consider  the following approximation
characteristics
\begin{equation}\label{velzag}
{\cal E}\left(C^{\psi}_{\beta,p};
Z_{n-1}^{s}\right)_{C}=\mathop{\sup}\limits_{f\in
C^{\psi}_{\beta,p}}\|f(\cdot)-Z^{s}_{n-1}(f;\cdot)\|_{C}, \ 1\leq
p\leq\infty, \ \beta\in
    \mathbb{R},
\end{equation}
and solve the problem of establishing the order of decreasing to
zero as $n\rightarrow\infty$ of the mentioned  quantities with
respect to relations between parameters $\psi$, $\beta$, $p$ and
$s$. It is clear that we can make conclusion about the approximation
ability of a linear polynomial approximation method (including Fejer
$\sigma_{n-1}$ and Zygmund  $Z_{n-1}^{s}$ methods) on the class
$C^{\psi}_{\beta,p}$, after comparison the rate of decreasing of the
exact upper bounds of uniform deviations of trigonometric sums,
which are generated  by this method, on the set $C^{\psi}_{\beta,p}$
with the rate of decreasing of the best uniform approximations of
the class $C^{\psi}_{\beta,p}$ by trigonometric polynomials
$t_{n-1}$ of order not higher than $n-1$, namely the quantities of
the form:

\begin{equation}\label{nnabl}
{E}_{n}(C^{\psi}_{\beta,p})_{C}=\sup\limits_{f\in
C^{\psi}_{\beta,p}}\inf\limits_{t_{n-1}}\|f(\cdot)-t_{n-1}(\cdot)\|_{C},
\ \ 1\leq p\leq\infty,
\end{equation}
where $\mathcal{T}_{2n-1}$  is the subspace of trigonometric
polynomials $t_{n-1}$ of order  $n-1$ with real coefficients. In
this case, since always the following estimate holds
\begin{equation}\label{ocnn}
{ E}_{n}\Big(C^{\psi}_{\beta,p}\Big)_{C}\leq {\cal
E}\left(C^{\psi}_{\beta,p}; Z_{n-1}^{s}\right)_{C}, \ n\in
\mathbb{N},
\end{equation}
it is important to  know under which restrictions on the parameters
$\psi, s, \beta$ and $p$ the following equality takes place
\begin{equation}\label{prnnabl}
{ E}_{n}\Big(C^{\psi}_{\beta,p}\Big)_{C}\asymp {\cal
E}\left(C^{\psi}_{\beta,p}; Z_{n-1}^{s}\right)_{C}.
\end{equation}
The notation $A(n)\asymp B(n)$ means, that  $A(n)=O(B(n))$  and at
the same time
 \linebreak$B(n)=O(A(n))$,  where by the notation
$A(n)=O(B(n))$  we mean, that there exists a constant $K>0$ such
that the inequality $A(n)\leq K(B(n))$ holds.

 In the work
\cite{Zi1} A. Zygmund introduced trigonometric polynomials of the
form (\ref{sz})  and found exact order estimates of the quantities
${\cal E}\left(W^{r}_{\infty};Z_{n-1}^{s}\right)_{C}$  at
$r\in\mathbb{N}$. B. Nagy investigated \cite{N}  the quantities
${\cal E}\left(W^{r}_{\beta,\infty};Z_{n-1}^{s}\right)_{C}$  at
$r>0$, $\beta\in \mathbb{Z}$,   and for $s\leq r$  he established
the asymptotic equality, and for $s> r$ he found order estimates.
Later S.A. Telyakovsky \cite{Telyakovskiy63} obtained asymptotically
exact equalities for the quantities ${\cal
E}\left(W^{r}_{\beta,\infty};Z_{n-1}^{s}\right)_{C}$ for $r>0$ and
$\beta\in \mathbb{R}$   for $n\rightarrow\infty$. On the Weyl-Nagy
classes, the exact order estimates of the quantities ${\cal
E}\left(W^{r}_{\beta,p};Z_{n-1}^{s}\right)_{C}$  for $1< p<\infty$
and $r>1/p$  and for $p=1$  and $r\geq1$, $\beta\in \mathbb{R}$
 are found in the work
    \cite{Kostuch1}.

     Concerning the Fejer sums  $\sigma_{n-1}(f;t)$  it should be noticed that the order estimates of quantities  ${\cal
E}\left(W^{r}_{\beta,\infty};\sigma_{n-1}\right)_{C}$, $r>0$  for
$\beta\in \mathbb{Z}$   were found by S.M. Nikol’skii  \cite{Nik};
 for the quantities ${\cal
E}\left(W^{r}_{r,p};\sigma_{n-1}\right)_{C}$ for $1<p\leq\infty$ and
$r>1/p$,  and also for $p=1$  and $r\geq1$  were found  by V.M.
Tikhomirov \cite{Tuxomurov}  and  by A.I. Kamzolov \cite{Kam2}.

 Approximation properties of Zygmund sums on the classes of  $(\psi,\beta)$-differentiable functions were studied in the works \cite{Bu}, \cite{SG}, \cite{Serduk1},
  (see., also,  \cite{S1}). Particularly in the work of
D.M. Bushev \cite{Bu} the asymptotic equalities for the quantities
${\cal E}(C^{\psi}_{\beta,\infty};Z^{s}_{n-1})_{C}$ were established
for some quite natural constraints on $\psi$ and $s$ as
$n\rightarrow\infty$. In the case, when the  series
$\sum_{k=1}^{\infty}\psi^{2}(k)$ is convergent, the exact values of
the quantities ${\cal E}\left(C^{\psi}_{\beta,2};Z_{n-1}^{s}
\right)_{C}$  were established in the work of A.S. Serdyuk and
I.V.Sokolenko \cite{Serduk1}.

In the work \cite{SG}  the authors found the exact order estimatites
of uniform approximations by Zygmund sums $Z^{s}_{n-1}$ on the
classes $C^{\psi}_{\beta,p}$, $1<p<\infty$, when
$\psi\in\Theta_{p}$, and $\Theta_{p}$, $1< p<\infty$, is the set of
non-increasing functions $\psi(t)$, for which there exists
$\alpha>1/p$ such that the function $t^{\alpha}\psi(t)$ almost
decreases, and $\psi(t)t^{s+1/p-\varepsilon}$ increases by
$[1,\infty)$ for some $\varepsilon>0$.

 Concerning the estimates of the best uniform approximations of functional compacts, it should be noticed the following. For the Weyl-Nagy classes
$W^{r}_{\beta,p}$, $r>1/p$, $\beta\in \mathbb{R}$, $1\leq
p\leq\infty$, the exact order estimates of the best approximations
${ E}_{n}\Big(W^{r}_{\beta,p}\Big)_{C}$ are known (see, e.g.,
\cite{Teml}). Moreover, for $p=\infty$ the exact values of the
quantities ${ E}_{n}\Big(W^{r}_{\beta,\infty}\Big)_{C}$ for all
$r>0$, $\beta\in \mathbb{R}$ and $n\in \mathbb{N}$  are known (see.
\cite{Dziadyk}).

The order estimates of the best approximations of the classes
$C^{\psi}_{\beta,p}$ under certain restrictions on $\psi$, $\beta$
and $p$ were investigated in the works  \cite{UMG}, \cite{SerStep1},
\cite{SerStep}, \cite{S1}.  In some partial cases  (especially for
$p=\infty$) the exact or asymptotically exact values of the
quantities ${ E}_{n}\Big(C^{\psi}_{\beta,p}\Big)_{C}$ (are also
known (see. \cite{Pink}, \cite{Serdyuk1995}, \cite{Serdyuk1999},
\cite{Serdyuk_2002_zb}, \cite{Serdyuk_2005_7}, \cite{Serdyuk_2020},
\cite{S1}).

In this paper, we establish the exact order estimates of the
quantities of the form  (\ref{velzag}) for all $1\leq p<\infty$ and
$\beta\in \mathbb{R}$,  in case, when
    $\psi(t)t^{1/p}\in\mathfrak{M}_0$,  the product
   $\psi(k)k^{s+1/p}$ generally monotonically increases,  $\psi(k)k^{s+1/p-\varepsilon}$ almost increases (according to Bernstein)
    for some
   $\varepsilon>0$ and for $1< p<\infty$
\begin{equation}\label{um2}
 \sum_{k=n}^{\infty}\psi^{p'}(k)k^{p'-2}<\infty, \
 \frac{1}{p}+\frac{1}{p'}=1,
\end{equation}
 and for $p=1$
\begin{equation}\label{ump1}
 \sum_{k=n}^{\infty}\psi(k)<\infty.
\end{equation}
The conditions (\ref{um2}) and (\ref{ump1}) and the monotonic
decreasing to zero of the sequence  $\psi(k)$ ensure the inclusion
of $\Psi_{\beta}\in L_{p'}$, $1/p+1/p'=1$, $1\leq p<\infty$ (see,
e.g., Lemma 12.6.6 from \cite{Z2}, s. 193.)

In this paper it is also shown  that for some conditions, Zygmund
sums (and at $s=1$ also the Fejer sums) realize the orders of the
best uniform approximations on the classes $C^{\psi}_{\beta,p}$,
that is, the order estimate
 (\ref{prnnabl}) is true.  Previously, this property was
proved for Fourier sums \cite{UMG}, \cite{SerStep}, \cite{StepS},
\cite{Step}.

Let us formulate some necessary definitions.

A non-negative sequence $a=\big\{a_{k}\big\}_{k=1}^{\infty}$, $k\in
\mathbb{N}$, is said  to be generally  monotonically increasing
(and write $a\in GM^{+}$), if there exists a constant $A\geq1$, such
that for any natural $n_{1}$ and $n_{2}$ such that $n_{1}\leq n_{2}$
inequalities are held
\begin{equation}\label{belov}
a_{n_{1}}+\sum_{k=n_{1}}^{m-1}|a_{k}-a_{k+1}|\leq A a_{m}, \
m=\overline{n_{1}, n_{2}}.
\end{equation}
It is easy to see that if the positive  sequence
$a=\big\{a_{k}\big\}_{k=1}^{\infty}$ increases, starting from some
number, then it generally  monotonically increasing.

A non-negative sequence  $a=\big\{a_{k}\big\}_{k=1}^{\infty}$, $k\in
\mathbb{N}$  is said  to be  almost increasing (according to
Bernstein) if there exists a constant  $K$, such that for all,
 $n_{1}\leq n_{2}$
\begin{equation}\label{mzr}
a_{n_{1}}\leq K a_{n_{2}}.
\end{equation}
In this case, if for the sequence
$a=\big\{a_{k}\big\}_{k=1}^{\infty}$ there exists a constant
$\varepsilon>0$, such that $\big\{a_{k}k^{-\varepsilon}\big\}$
almost increases, then we write $a\in GA^{+}$. It is clear that if
the sequence $a$ belongs to $GM^{+}$,  then it is almost increasing
according to Bernstein.

 Let us put further at   $\delta>0$  $g_{\delta}(t):=\psi(t)t^{\delta}$,
 $t\in[1,\infty)$.

\vskip 3mm

\section{Order estimates of the  approximations  by Zygmund sums on the classes of convolutions}

\bf Theorem 1.  \it{Let   $s>0$, $1\leq p<\infty$, $g_{1/p}\in
{\mathfrak M}_{0}$,  $g_{s+1/p}\in GM^{+}\cap GA^{+}$, $\beta\in
\mathbb{R}$ and $n\in \mathbb{N}$. In the case $1<p<\infty$, if the
condition (\ref{um2}) holds and the following inequality holds
\begin{equation}\label{um3}
\inf\limits_{t\geq1}\alpha(g_{1/p};t)>\frac{p'}{2},
\end{equation}
then the following order estimates take place
\begin{equation}\label{t5}
{E}_{n}\left(C^{\psi}_{\beta,p}\right)_{C}\asymp{\cal
E}\left(C^{\psi}_{\beta,p};
Z_{n-1}^{s}\right)_{C}\asymp\left(\sum_{k=n}^{\infty}\psi^{p'}(k)k^{p'-2}\right)^{1/p'},
\ \frac{1}{p}+\frac{1}{p'}=1;
\end{equation}
in the case $p=1$, if   the condition (\ref{ump1}) holds and the
following inequality holds
\begin{equation}\label{um3p1}
\inf\limits_{t\geq1}\alpha(g_{1};t)>1,
\end{equation}
then the following order estimates take place
\begin{equation}\label{ocp1}
{E}_{n}\left(C^{\psi}_{\beta,1}\right)_{C}\asymp{\cal
E}\left(C^{\psi}_{\beta,1}; Z_{n-1}^{s}\right)_{C}\asymp
{\left\{{\begin{array}{l l}
\sum\limits_{k=n}^{\infty}\psi(k), & \cos \frac{\beta\pi}{2}\neq 0;  \\
\psi(n)n, &\cos \frac{\beta\pi}{2}= 0.
\end{array}} \right.}
\end{equation}

\textbf{Proof.}  \rm  Since the operator $Z^{s}_{n-1}:
f(t)\rightarrow Z^{s}_{n-1}(f,t)$ is linear polynomial operator,
which is invariant under the shift, i.e.
$$
Z^{s}_{n-1}(f_{h}, t)= Z^{s}_{n-1}(f, t+h), \ f_{h}(t)=f(t+h), \
h\in \mathbb{R},
$$
and norm  in $C$ and classes $C^{\psi}_{\beta,p}$ also are invariant
under the shift, that is
$$
\|f_{h}(t)\|_{C}=\|f(t)\|_{C}; \ f(t)\in
C^{\psi}_{\beta,p}\Rightarrow f_{h}(t)\in C^{\psi}_{\beta,p},
$$
then
\begin{equation}\label{nornol}
{\cal E}\left(C^{\psi}_{\beta,p};
Z_{n-1}^{s}\right)_{C}=\mathop{\sup}\limits_{f\in C^{\psi}_{\beta,p}
}|f(0)-Z_{n-1}^{s}(f;0)|.
\end{equation}
By virtue  (\ref{zgo}) and (\ref{sz}) for any   function
 $f\in C^{\psi}_{\beta,p}$, $1\leq p<\infty$, $\beta\in
\mathbb{R}$, $s>0$ the following equality holds
\begin{equation}\label{predst}
f(0)-Z_{n-1}^{s}(f;0)=\frac{1}{\pi}\!\!\int\limits_{-\pi}^{\pi}\!\!\!
\left(\frac{1}{n^{s}}\sum_{k=1}^{n-1}\psi(k)k^{s}\cos\left(kt+\frac{\beta\pi}{2}\right)\!+\!\Psi_{-\beta,n
    }(t)\right)\!\varphi(t)dt,
\end{equation}
where $\Psi_{-\beta,n
    }(t)=\sum_{k=n}^{\infty}\psi(k)\cos\left(kt+\frac{\beta\pi}{2}\right)$,
\ \  $\|\varphi\|_{p}\leq1$, $n\in \mathbb{N}$.

Relations (\ref{nornol}) and (\ref{predst}),
 H\"{o}lder's inequality and triangle inequality
imply that for  \linebreak$1\leq p<\infty$
\begin{equation}\label{111}
{\cal E}\left(C^{\psi}_{\beta,p};
Z_{n-1}^{s}\right)_{C}\leq\frac{1}{\pi
}\Bigg\|\frac{1}{n^{s}}\sum_{k=1}^{n-1}\psi(k)k^{s}\cos\left(kt+\frac{\beta\pi}{2}\right)+
\Psi_{-\beta,n
    }(t)\Bigg\|_{p'}\leq
$$
$$
\leq\frac{1}{\pi n^{s}}\Bigg\|
\sum_{k=1}^{n-1}\psi(k)k^{s}\cos\left(kt+\frac{\beta\pi}{2}\right)\Bigg\|_{p'}+\frac{1}{\pi}\big\|\Psi_{-\beta,n}(t)\big\|_{p'},
\ \frac{1}{p}+\frac{1}{p'}=1.
\end{equation}

Let us show that,  if  $g_{s+1/p}\in GM^{+}\cap GA^{+}$, where
$g_{s+1/p}=\big\{\psi(k)k^{s+1/p}\big\}_{k=1}^{\infty}$, then
\begin{equation}\label{v1a+}
  \Big\|
\sum_{k=1}^{n-1}\psi(k)k^{s}\cos\left(kt+\frac{\beta\pi}{2}\right)\Big\|_{p'}=O\big(\psi(n)n^{s+\frac{1}{p}}\big),
\ 1\leq p<\infty.
\end{equation}

Applying Abel transformation to the function
$\sum_{k=1}^{n-1}\psi(k)k^{s}\cos\left(kt+\frac{\beta\pi}{2}\right)$,
 we have
\begin{equation}\label{dkb}
\sum_{k=1}^{n-1}\psi(k)k^{s}\cos\Big(kt+\frac{\beta\pi}{2}\Big)=\sum_{k=1}^{n-2}\Big(\psi(k)k^{s}-
\psi(k+1)(k+1)^{s}\Big)D_{k,\beta}(t)+
$$
$$
+\psi(n-1)(n-1)^{s}D_{n-1,\beta}(t)-\frac{1}{2}\cos\frac{\beta\pi}{2},
\end{equation}
 where
$$
D_{k,\beta}(t):=\frac{1}{2}\cos\frac{\beta\pi}{2}+\sum_{\nu=1}^{k}\cos\Big(\nu
t-\frac{\beta\pi}{2}\Big).
$$
Then, in view of
$$
\|D_{k,\beta}(t)\|_{p'}=O(k^{1-\frac{1}{p'}})=O(k^{\frac{1}{p}}), \
1\leq p <\infty, \ k\in \mathbb{N}, \ \beta\in \mathbb{R}
$$
(see, e.g., \cite{UMG}), of (\ref{dkb}) we get
\begin{equation}\label{ddzn}
\Bigg\|\sum_{k=1}^{n-1}\psi(k)k^{s}\cos\left(kt+\frac{\beta\pi}{2}\right)\Bigg\|_{p'}=
$$
$$
=O(1)+O\Bigg(\sum_{k=1}^{n-2}
\big|\psi(k)k^{s}-\psi(k+1)(k+1)^{s}\big|k^{\frac{1}{p}}\Bigg)+
O\left(\psi(n-1)(n-1)^{s+\frac{1}{p}}\right).
\end{equation}

Since  $g_{s+1/p}\in GM^{+}$, then, by using the triangle
inequality,  inequality (\ref{belov}) and Lagrange theorem, we have
\begin{equation}\label{ddzn1}
\sum_{k=1}^{n-2}\!
\big|\psi(k)k^{s}-\psi(k+1)(k+1)^{s}\big|k^{\frac{1}{p}}\leq
$$
$$
\leq\sum_{k=1}^{n-2}
\big|\psi(k)k^{s+\frac{1}{p}}-\psi(k+1)(k+1)^{s+\frac{1}{p}}\big|+
\sum_{k=1}^{n-2}\big|\psi(k+1)(k+1)^{s+\frac{1}{p}}-\psi(k+1)(k+1)^{s}k^{\frac{1}{p}}\big|\leq
$$
$$
 \leq A\psi(n-1)(n-1)^{s+\frac{1}{p}}+
\frac{1}{p}\sum_{k=1}^{n-2}\psi(k+1)(k+1)^{s}k^{\frac{1}{p}-1}=
$$
$$
= A\psi(n-1)(n-1)^{s+\frac{1}{p}}+
\frac{1}{p}\sum_{k=1}^{n-2}\psi(k+1)(k+1)^{s+\frac{1}{p}-1}\big(1+\frac{1}{k}\big)^{\frac{1}{p'}}\leq
$$
$$
\leq A\psi(n-1)(n-1)^{s+\frac{1}{p}}+
2\sum_{k=2}^{n-1}\frac{\psi(k)k^{s+\frac{1}{p}}}{k}.
\end{equation}

 According to the condition $g_{s+1/p}\in GA^{+}$,
there exits  $\varepsilon>0$ such that the sequence
\linebreak$\big\{g_{s+1/p}(k)k^{-\varepsilon}\big\}=\big\{\psi(k)k^{s+1/p-\varepsilon}\big\}$
  almost increases, and hence
taking into account (\ref{mzr}), we obtain
\begin{equation}\label{ocsumy}
\sum\limits_{k=2}^{n-1}\frac{\psi(k)k^{s+1/p}}{k}=\sum\limits_{k=2}^{n-1}\frac{\psi(k)k^{s+1/p-\varepsilon}}{k^{1-\varepsilon}}\leq
$$
$$
\leq
K\psi(n-1)(n-1)^{s+1/p-\varepsilon}\sum_{k=2}^{n-1}\frac{1}{k^{1-\varepsilon}}<
K\psi(n-1)(n-1)^{s+1/p-\varepsilon}\int\limits_{1}^{n-1}\frac{dt}{t}<\frac{K}{\varepsilon}\psi(n-1)(n-1)^{s+1/p}.
\end{equation}
From (\ref{ddzn1}) and (\ref{ocsumy})  we get the following
inequality
\begin{equation}\label{ostner}
\big|\psi(k)k^{s}-\psi(k+1\!)(k+1)^{s}\big|k^{\frac{1}{p}}\leq\Big(A+\frac{2K}{\varepsilon}\Big)\psi(n-1)(n-1)^{s+1/p}.
\end{equation}
From (\ref{ddzn}) and (\ref{ostner}) we obtain an estimate
(\ref{v1a+}).

To estimate the norm $\|\Psi_{-\beta,n}(\cdot)\|_{p'}$ for
$1<p'<\infty$  we use the statement, which was  established in
\cite{SerStep}, and according to which in the case when
$\big\{a_{k}\big\}_{k=1}^{\infty}$ is
 the monotonically non-increasing sequence of positive
 numbers is such that
 $\sum_{k=1}^{\infty}a_{k}^{p'}k^{p'-2}<\infty$, then for an arbitrary $n\in \mathbb{N}$ and $\gamma\in\mathbb{R}$  the following estimate holds
\begin{equation}\label{nsumy}
\Big\|\sum_{k=n}^{\infty}a_{k}\cos\big(kx+\gamma\big)\Big\|_{p'}=O\Big(\sum_{k=n}^{\infty}a_{k}^{p'}k^{p'-2}+a_{n}^{p'}n^{p'-1}\Big)^{1/p'}.
\end{equation}
 Putting in
(\ref{nsumy}) $a_{k}=\psi(k)$, $\gamma=\frac{\beta\pi}{2}$
 we obtain that for  $1<p<\infty$,
$\beta\in\mathbb{R}$ and $n\in \mathbb{N}$
\begin{equation}\label{nb}
\|\Psi_{-\beta,n
    }(\cdot)\|_{p'}=O\Big(\sum_{k=n}^{\infty}\psi^{p'}(k)k^{p'-2}+\psi^{p'}(n)n^{p'-1}\Big)^{1/p'}.
\end{equation}

 Then, using Lemma 3 of
\cite{SerStep}, we conclude that for $1<p'<\infty$,  $n\in
\mathbb{N}$ under condition (\ref{um2}) and imbedding $g_{1/p}\in
{\mathfrak M_{0}}$ the following estimate holds
\begin{equation}\label{poc0}
\psi^{p'}(n)n^{p'-1}=O\Bigg(\sum_{k=n}^{\infty}\psi^{p'}(k)k^{p'-2}\Bigg).
\end{equation}
 According to the conditions of Theorem 1 we have that
$g_{1/p}\in {\mathfrak M_{0}}$,  so  taking into account
(\ref{poc0}), from (\ref{nb}),  we obtain
\begin{equation}\label{nb1}
\|\Psi_{-\beta,n
    }(\cdot)\|_{p'}=O\Bigg(\sum_{k=n}^{\infty}\psi^{p'}(k)k^{p'-2}\Bigg)^{1/p'}, \ 1<
    p'<\infty, \ \beta\in\mathbb{R}, \ n\in \mathbb{N}.
\end{equation}

 Combining (\ref{111}), (\ref{v1a+}) and (\ref{nb1})  in the case when  $g_{1/p}\in {\mathfrak M_{0}}$, and
$g_{s+1/p}\in GM^{+}\cap GA^{+}$,
  we arrive at the estimate
\begin{equation}\label{oczv}
{\cal E}\left(C^{\psi}_{\beta,p}; Z_{n-1}^{s}\right)_{C}=
O\Bigg(\sum_{k=n}^{\infty}\psi^{p'}(k)k^{p'-2}\Bigg)^{1/p'}, \ \
1<p<\infty, \ \frac{1}{p}+\frac{1}{p'}=1.
\end{equation}

 As follows from Corollary 1 and 2 from  \cite{SerStep}  for
$1<p<\infty$, $1/p+1/p'=1$, $n\in \mathbb{N}$  and $\beta\in
\mathbb{R}$, under  conditions  (\ref{um2}) and (\ref{um3})  and
imbedding $g_{1/p}\in {\mathfrak M_{0}}$ for ${
E}_{n}\Big(C^{\psi}_{\beta,p}\Big)_{C}$
 we arrive at
the following  order estimates
\begin{equation}\label{ocznnn}
{ E}_{n}\Big(C^{\psi}_{\beta,p}\Big)_{C}\asymp
\Big(\sum_{k=n}^{\infty}\psi^{p'}(k)k^{p'-2}\Big)^{1/p'}.
\end{equation}

\rm  Therefore, by virtue of inequality (\ref{ocnn})  and relations
(\ref{oczv}) and (\ref{ocznnn})   we obtain   order equality
(\ref{t5}).

 Further, let us consider the case   $p=1$.
 Let us establish the estimate of the norm
\linebreak$\|\Psi_{-\beta,n
    }(\cdot)\|_{p'}=\|\Psi_{-\beta,n
    }(\cdot)\|_{\infty}$.

 It is obvious that  for any  $\beta\in \mathbb{R}$ the following inequality holds

\begin{equation}\label{nbkos}\|\Psi_{-\beta,n
    }(\cdot)\|_{\infty}=
    \bigg\|\sum_{k=n}^{\infty}\psi(k)\cos\left(kt+\frac{\beta\pi}{2}\right)\bigg\|_{\infty}\leq \sum_{k=n}^{\infty}\psi(k).
\end{equation}
      If $\beta=2k+1$, $k\in \mathbb{Z}$,  then following estimate takes place
\begin{equation}\label{nbsin}
\|\Psi_{-\beta,n
    }(\cdot)\|_{\infty}=
    \bigg\|\sum_{k=n}^{\infty}\psi(k)\sin
    kt\bigg\|_{\infty}\leq(\pi+2)\psi(n)n
\end{equation}
(see, e.g., relation
    (82)  from \cite{Step}).

 According to Lemma 3   from \cite{Step}, if  $g_{1}\in \mathfrak{M}_0$,  where $g_{1}=\big\{\psi(k)k\big\}_{k=1}^{\infty}$ and the  condition (\ref{ump1})
holds,
 then the
following estimates are true
\begin{equation}\label{poc2}
\psi(n)n=O\bigg(\sum_{k=n}^{\infty}\psi(k)\bigg).
\end{equation}

If  $g_{1}\in \mathfrak{M}_0$ and   the conditions (\ref{ump1})
hold, then combining (\ref{111}), (\ref{v1a+}), (\ref{nbkos}) --
(\ref{poc2}),
 we obtain the following estimates
\begin{equation}\label{oczv1}
{\cal E}\left(C^{\psi}_{\beta,1}; Z_{n-1}^{s}\right)_{C}={\left\{
{\begin{array}{l l}
O\Big(\sum\limits_{k=n}^{\infty}\psi(k)\Big), & \cos \frac{\beta\pi}{2}\neq 0; \\
O\big(\psi(n)n\big), &\cos \frac{\beta\pi}{2}= 0.
\end{array}} \right.}
\end{equation}

To estimate the quantity  ${\cal E}\left(C^{\psi}_{\beta,1};
Z_{n-1}^{s}\right)_{C}$    from below, we use  Theorems 3  and 4
  from \cite{Step},  according to which, if $g_{1}\in \mathfrak{M}_0$ and the conditions
 (\ref{ump1})  and (\ref{um3p1}) are true, then   for $n\in \mathbb{N}$  and $\beta\in
\mathbb{R}$   the following the order equalities take place
\begin{equation}\label{oczn1}
{ E}_{n}\Big(C^{\psi}_{\beta,1}\Big)_{C}\asymp{\left\{
{\begin{array}{l l}
\sum\limits_{k=n}^{\infty}\psi(k), & \cos \frac{\beta\pi}{2}\neq 0; \\
\psi(n)n, &\cos \frac{\beta\pi}{2}= 0.
\end{array}} \right.}
\end{equation}
The estimate  (\ref{ocp1})  follows from the inequality
(\ref{ocnn}), estimates (\ref{oczv1}) and (\ref{oczn1}).  Theorem 1
is proved.

Assume that the conditions of Theorem 1  take place,   moreover,
more stronger imbedding holds   $g_{1/p}\in \mathfrak{M}_{C}$.  As
it follows from Lemma 3  from \cite{SerStep}  if $g_{1/p}\in
\mathfrak{M}_{C}$ and  the condition  (\ref{um2})  holds, then  for
$1<p<\infty$ the following estimates take place
\begin{equation}\label{pocc}
\sum_{k=n}^{\infty}\psi^{p'}(k)k^{p'-2}\asymp\psi^{p'}(n)n^{p'-1}.
\end{equation}
In addition, as  it was shown in \cite{Step}, Lemma 3], if
$g_{1}\in \mathfrak{M}_{C}$ and the condition (\ref{ump1})  holds,
then
 the following order estimates are true
\begin{equation}\label{poc21}
\sum_{k=n}^{\infty}\psi(k)\asymp\psi(n)n.
\end{equation}
Formulas (\ref{pocc})  and (\ref{poc21}), and Theorem 1 allow us to
write the following statement.

\bf Theorem 2.  \it{Let Let $s>0$, $1\leq p<\infty$, $g_{1/p}\in
\mathfrak{M}_{C}$, $g_{s+1/p}\in GM^{+}\cap GA^{+}$, $\beta\in
\mathbb{R}$  and $n\in \mathbb{N}$.

In the case $1<p<\infty$, if  the conditions  (\ref{um2}) and
(\ref{um3}) hold, then the following order estimates take place
\begin{equation}\label{t5c}
{E}_{n}(C^{\psi}_{\beta,p})_{C}\asymp{\cal
E}\left(C^{\psi}_{\beta,p};
Z_{n-1}^{s}\right)_{C}\asymp\psi(n)n^{1/p},
\end{equation}
 and in the case  $p=1$  if  the condition  (\ref{ump1})  and (\ref{um3p1}) hold, then
 the following order estimates take place
\begin{equation}\label{ocp1c}
{E}_{n}(C^{\psi}_{\beta,1})_{C}\asymp{\cal
E}\left(C^{\psi}_{\beta,1}; Z_{n-1}^{s}\right)_{C}\asymp\psi(n)n.
\end{equation}

\textbf{Proof.} \rm Order estimates  (\ref{t5c}) were established in
\cite{SG}.

  Note that when  $1< p<\infty$, $g_{1/p}\in \mathfrak{M}_{0}$  and
\begin{equation}\label{alfan}
\lim\limits_{t\rightarrow\infty}\alpha\big(g_{1/p};t\big)=\infty,
\end{equation}
then  the order estimates  (\ref{t5c}) do not take place,   since in
this case
$$
\psi(n)n^{\frac{1}{p}}=o\Big(\big(\sum_{k=n}^{\infty}\psi^{p'}(k)k^{p'-2}\big)^{1/p'}\Big),
\ n\rightarrow\infty
$$
(see, Lemma from \cite{SerStep}).

 Similarly, when $p=1$, $g_{1/p}=g_{1}\in
\mathfrak{M}_{0}$  and
\begin{equation}\label{alfan1}
\lim\limits_{t\rightarrow\infty}\alpha\big(g_{1};t\big)=\infty,
\end{equation}
then as follows from Lemma 3 \cite{Step}
$$
\psi(n)n=o\Big(\sum_{k=n}^{\infty}\psi(k)\Big),
$$
in this case,  for    $\beta$ such that
 $\cos \frac{\beta\pi}{2}\neq 0$
  order estimates (\ref{ocp1c}) do not take place.

 As example of the function  $\psi(t)$, for which the conditions of  Theorem 1
   and the equalities (\ref{alfan})  and
(\ref{alfan1}) take place, we can use the function
\begin{equation}\label{funk}
\psi(t)=t^{-1/p}\ln^{-\gamma}(t+K), \ \gamma>{\left\{
{\begin{array}{l l} \frac{1}{p'}, &1<p<\infty; \\  1, & p=1,
\end{array}} \right.} \
K>{\left\{
{\begin{array}{l l}e^{\gamma p'/2}, & 1<p<\infty; \\
e^{\gamma}, & p=1,
\end{array}} \right.}
\end{equation}
(see \cite{SerStep}, \cite{Step}). Let us write  the order estimates
for the quantities
 ${E}_{n}\Big(C^{\psi}_{\beta,p}\Big)_{C}$ and \linebreak${\cal
E}\left(C^{\psi}_{\beta,p}; Z_{n-1}^{s}\right)_{C}$  in the case,
when  $\psi(t)$  has the form
 (\ref{funk}).

\bf Theorem 3.  \it{ Let $\psi(t)=t^{-1/p}\ln^{-\gamma}(t+K)$,
$\beta\in \mathbb{R}$ and $n\in \mathbb{N}$. If $1<p<\infty$,
$\gamma>1/p'$, $K>e^{\gamma p'/2}$, $1/p+1/p'=1$, then
\begin{equation}\label{pof}
{E}_{n}(C^{\psi}_{\beta,p})_{C}\asymp{\cal
E}\left(C^{\psi}_{\beta,p};
Z_{n-1}^{s}\right)_{C}\asymp\psi(n)n^{1/p}\ln^{1/p'}n, \ n\geq2;
\end{equation}
if  $p=1$, $\gamma>1$, $K> e^{\gamma}$, then
\begin{equation}\label{ocp1f}
{E}_{n}(C^{\psi}_{\beta,1})_{C}\asymp{\cal
E}\left(C^{\psi}_{\beta,1}; Z_{n-1}^{s}\right)_{C}\asymp{\left\{
{\begin{array}{l l}
\psi(n)n\ln n, & \cos \frac{\beta\pi}{2}\neq 0, \ \  n\geq2; \\
\psi(n)n, &\cos \frac{\beta\pi}{2}= 0.
\end{array}} \right.}
\end{equation}
\rm

  We show that for the indicated function  $\psi$ of the form  (\ref{funk}) all conditions of the Theorem 1 are true.  Indeed,  for  $1<p<\infty$,
  $\gamma>1/p'$,
$K>e^{\gamma p'/2}$  we have
$$\sum\limits_{k=n}^{\infty}\psi^{p'}(k)k^{p'-2}=\sum_{k=n}^{\infty}\frac{1}{k\ln^{\gamma
p'}(k+K)}<\infty,$$
$$\alpha\big(g_{1/p};t\big)=\frac{(t+K)\ln(t+K)}{\gamma t}>\frac{\ln(t+e^{\gamma p'/2})}{\gamma
},$$ and hence,
$\lim\limits_{t\rightarrow\infty}\alpha\big(g_{1/p};t\big)=\infty$
 and $\alpha\big(g_{1/p};t\big)>\frac{p'}{2}.$

 For  $p=1$, $\gamma>1$, $K\geq e^{\gamma}$,  we have
$\sum\limits_{k=n}^{\infty}\psi(k)\leq\sum\limits_{k=n}^{\infty}\frac{1}{k\ln^{\gamma
}(k+e^{\gamma})}<\infty$,
$$\alpha\big(g_{1};t\big)>\frac{\ln(t+e^{\gamma})}{\gamma
},$$ and hence,
$\lim\limits_{t\rightarrow\infty}\alpha\big(g_{1};t\big)=\infty$
 і $\alpha\big(g_{1};t\big)>1.$

It is obvious that for any   $s>0$ and $1\leq p<\infty$ the
functions $g_{s+1/p}(t)=t^{s}\ln^{-\gamma}(t+K)$  increase
monotonically, starting from some point  $t_{0}$. Therefore, it is
not difficult to be convinced that the sequence $g_{s+1/p}(k)$
belongs to the set $GM^{+}\cap GA^{+}$

Therefore, the  function $\psi$ of the form (\ref{funk})
 satisfies the conditions of Theorem 1.

Further,  using the formula (79)  from \cite{SerStep}, obtain
\begin{equation}\label{ocfunkd}
\Big(\sum_{k=n}^{\infty}\psi^{p'}(k)k^{p'-2}\Big)^{1/p'}\asymp\Big(\int\limits_{n}^{\infty}\psi^{p'}(t)t^{p'-2}dt\Big)^{1/p'}=
\Big(\int\limits_{n}^{\infty}\frac{dt}{t\ln^{\gamma
p'}(t+K)}\Big)^{1/p'}
$$
$$
\asymp\ln^{1/p'-\gamma}n=\psi(n)n^{1/p}\ln^{1/p'}n\frac{\ln^{-\gamma}n}{\ln^{-\gamma}(n+K)}\asymp\psi(n)n^{1/p}\ln^{1/p'}n,
\ n\geq2.
\end{equation}
 Then formula (\ref{pof})  follows from the estimate (\ref{t5})  and
 relations (\ref{ocfunkd}).

Similarly, by virtue of the inequality  (87)  from \cite{Step}
 we get
\begin{equation}\label{ocfunkd1}
\sum_{k=n}^{\infty}\psi(k)\asymp\int\limits_{n}^{\infty}\psi(t)dt=
\int\limits_{n}^{\infty}\frac{dt}{t\ln^{\gamma}(t+K)}
$$
$$
\asymp\ln^{1-\gamma}n\asymp\psi(n)n\ln n, \ n>2.
\end{equation}

 Formula  (\ref{ocp1f})  follows from the estimates  (\ref{ocp1})
 and relations (\ref{ocfunkd1}), in the case where $\beta$ is
such that $\cos \frac{\beta\pi}{2}\neq 0$.  By this Theorem 3 is
proved.

As it was already mentioned, for  $s=1$   the sums Zygmund
$Z_{n-1}^{s}$ coincide with the known Fejer sums $\sigma_{n-1}$.
Therefore, Theorem 1 and 2 imply the following statements.

\bf Proposition 1.  \it{Let  $1\leq p<\infty$, $g_{1/p}\in
\mathfrak{M}_{0}$,  $g_{1+1/p}\in GM^{+}\cap GA^{+}$, $\beta\in
\mathbb{R}$ and $n\in \mathbb{N}$. In the case $1<p<\infty$,  if the
conditions (\ref{um2}) and (\ref{um3}) hold, then the following
order estimates take place
\begin{equation}\label{t5sf}
{E}_{n}(C^{\psi}_{\beta,p})_{C}\asymp{\cal
E}\left(C^{\psi}_{\beta,p};
\sigma_{n-1}\right)_{C}\asymp\big(\sum_{k=n}^{\infty}\psi^{p'}(k)k^{p'-2}\big)^{1/p'};
\end{equation}
in the case $p=1$, if the conditions (\ref{ump1}) and (\ref{um3p1})
hold, then the following order equlaities take place
\begin{equation}\label{ocp1s}
{E}_{n}(C^{\psi}_{\beta,1})_{C}\asymp{\cal
E}\left(C^{\psi}_{\beta,1}; \sigma_{n-1}\right)_{C}\asymp
{\left\{{\begin{array}{l l}
\sum\limits_{k=n}^{\infty}\psi(k), & \cos \frac{\beta\pi}{2}\neq 0,  \\
\psi(n)n, &\cos \frac{\beta\pi}{2}= 0.
\end{array}} \right.}
\end{equation}

\bf Proposition 2.  \it{ Let $1\leq p<\infty$, $g_{1/p}\in
\mathfrak{M}_{C}$,  $g_{1+1/p}\in GM^{+}\cap GA^{+}$, $\beta\in
\mathbb{R}$ and $n\in \mathbb{N}$. In the case $1<p<\infty$, if  the
conditions  (\ref{um2}) and (\ref{um3}) hold, then the following
order estimates take place
\begin{equation}\label{t5sf1}
{E}_{n}(C^{\psi}_{\beta,p})_{C}\asymp{\cal
E}\left(C^{\psi}_{\beta,p};
\sigma_{n-1}\right)_{C}\asymp\psi(n)n^{1/p};
\end{equation}
in the case  $p=1$, if  the conditions  (\ref{ump1}) and
(\ref{um3p1}) hold, then the following order estimates take place
\begin{equation}\label{ocp1s}
{E}_{n}(C^{\psi}_{\beta,1})_{C}\asymp{\cal
E}\left(C^{\psi}_{\beta,1}; \sigma_{n-1}\right)_{C}\asymp\psi(n)n.
\end{equation}

\bf {References} \rm

\end{document}